\documentclass{amsart}
\usepackage{amsfonts,amssymb}
\usepackage{graphicx,psfrag}

\begin{document}


\newcommand{\RR}{\mathbb{R}} 
\newcommand{\ZZ}{\mathbb{Z}}
\newcommand{\QQ}{\mathbb{Q}}
\newcommand{\CC}{\mathbb{C}}
\newcommand{\EE}{\mathbb{E}}
\newcommand{\NN}{\mathbb{N}}
\newcommand{\LL}{\mathbb{L}}
\newcommand{\cA}{\mathcal{A}}
\newcommand{\cC}{\mathcal{C}}
\newcommand{\cF}{\mathcal{F}}
\newcommand{\cS}{\mathcal{S}}
\newcommand{\cM}{\mathcal{M}}
\newcommand{\cR}{\mathcal{R}}
\newcommand{\cD}{\mathcal{D}}
\newcommand{\cO}{\mathcal{O}}
\newcommand{\cT}{\mathcal{T}}
\newcommand{\cE}{\mathcal{E}}
\newcommand{\cG}{\mathcal{G}}
\newcommand{\cZ}{\mathcal{Z}}
\newcommand{\ux}{\mathbf{x}}
\newcommand{\uy}{\mathbf{y}}
\newcommand{\uz}{\mathbf{z}}
\newcommand{\uw}{\mathbf{w}}
\def\bbone{{\mathchoice {\rm 1\mskip-4mu l} {\rm 1\mskip-4mu l}
{\rm 1\mskip-4.5mu l} {\rm 1\mskip-5mu l}}}

\author{Abdelmalek Abdesselam}
\address{Abdelmalek Abdesselam, Department of Mathematics,
P. O. Box 400137,
University of Virginia,
Charlottesville, VA 22904-4137, USA}
\email{malek@virginia.edu}

\title{Towards Three-Dimensional Conformal Probability}

\begin{abstract}
In this outline of a program, based on rigorous renormalization group theory,
we introduce new definitions which allow one to formulate precise mathematical
conjectures related to conformal invariance as studied by physicists in the
area known as higher-dimensional conformal
bootstrap which has developed at a breathtaking pace over the last few years.
We also explore a second theme, intimately tied to conformal invariance for random distributions, which can be construed as a search for
very general first and second-quantized Kolmogorov-Chentsov Theorems. First-quantized refers to regularity of sample paths.
Second-quantized refers to regularity of generalized functionals or Hida distributions and relates to the operator product expansion.
We formulate this program in both the Archimedean and $p$-adic situations. 
Indeed, the study of conformal field theory and its connections with probability provides a golden opportunity where $p$-adic analysis can lead the way towards a better understanding of open problems in the Archimedean setting.
Finally, we present a summary of progress made on a $p$-adic hierarchical model and point out possible connections to number theory.
Parts of this article were presented in author's talk at the 6th International Conference on $p$-adic Mathematical Physics and its Applications, Mexico 2017.
\end{abstract}

\maketitle

\section{General context}\label{context}

The typical calculus sequence in North American universities goes as follows: Calculus I (differentiation in one variable), Calculus II
(integration in one variable), Calculus III (differentiation and integration in finitely many variables). As is well known,
the programmatic content of this sequence was worked out at the heuristic level by Newton, Leibniz and others in the 17th century.
Yet it almost took two more centuries in order for it to become a fully rigorous part of mathematics. 
One may argue that this sequence continues well into the graduate curriculum, and far beyond, into
the unknown. Calculus IV (differentiation in infinitely many variables) relates to the calculus of variations and the
notions of Gateau and Fr\'echet differentiability. It became well understood in the beginning of the 20th century.
By contrast, Calculus V (integration in infinitely many variables or functional integration) is at a less advanced stage.
It started with the construction of Wiener measure for Brownian motion and most of modern probability theory or stochastic calculus
can be seen as an effort to develop
Calculus V. A deep theory for functional integration has been invented by Kenneth G. Wilson, namely,
{\em renormalization group} (RG) theory. It is largely heuristic. However, as emphasized in~\cite{BrydgesPrague} it is a {\em systematic}
calculus with a wide range of applications. Indeed, a search on Google Scholar with the exact phrase ``renormalization group'' returns
more than 206,000 articles! Most of this literature is not mathematically rigorous and belongs to physics, but there is also a substantial
body of work where the RG is implemented in a mathematically rigorous way. The contents of this program
fall within this research area known as rigorous RG theory or modern constructive quantum field theory.
Rigorous RG theory started with the results~\cite{BleherS1,BleherS2} on hierarchical models and it was shown to work in
the Euclidean case in~\cite{GawedzkiKgr3}. 
Recent results in the field are: 1) Falco's work on the Coulomb gas in 2D~\cite{Falco1,Falco2}; 2)
the series of 
articles~\cite{BauerschmidtBS1,BauerschmidtBS2,BauerschmidtBS3,BauerschmidtBS4,BauerschmidtBS5,BrydgesS1,BrydgesS2,BrydgesS3,BrydgesS4,SladeT}
by Bauerschmidt, Brydges, Slade and Tomberg on the $N$-component $\phi^4$
model in four dimensions, including $N=0$ or the weakly self-avoiding walk; 3) the work by Chandra,
Guadagni and the author~\cite{AbdesselamCG0,AbdesselamCG1,AbdesselamCG2,Chandra}
on the $p$-adic fractional $\phi_{3}^{4}$ model; and more recently 4) the work of Lohmann, Slade and Wallace on the Archimedean fractional $\phi_{3}^{4}$ model~\cite{Slade,LohmannSW}.
The items 3) and 4) are particularly challenging
since the regime studied is below the upper critical dimension and involves a nontrivial RG fixed point corresponding to a highly {\em non-Gaussian} scaling limit (the word ``highly'' will be explained later).

{\em Conformal probability} (CP) is a very active area in analysis. Its main objects of study are scaling limits of models from
statistical mechanics. Such limits often exhibit a rich and nontrivial symmetry, namely, {\em conformal invariance} (CI).
The simplest example is provided
by the time-inversion invariance of Brownian motion $t B_{1/t}\stackrel{d}{=}B_t$ proved by Paul L\'evy in~\cite{Levy}.
Much deeper is the analogous invariance under conformal transformation in the 2D space of ``time parameters'' for
the Ising model critical scaling limit~\cite{CamiaGN,ChelkakHI} which represents the culmination of major effort led by Schramm, Smirnov
and others. Define the lattice couplings $J_{\ux,\uy}$ for $\ux,\uy\in\ZZ^2$ by the indicator function of $\ux,\uy$ being nearest neighbors.
Given a single site probability measure $\rho_{\rm ss}=\frac{1}{2}(\delta_{-1}+\delta_{1})$ on $\{-1,1\}$, a nonnegative parameter $\beta$,
and a finite set or volume $\Lambda\subset\ZZ^2$, one can define the probability measure
$\nu_{\beta,\Lambda}$ on $\{-1,1\}^{\ZZ^2}$ by
\[
{\rm d}\nu_{\beta,\Lambda}(\sigma)=\frac{1}{\cZ_{\beta,\Lambda}}
\exp\left(\frac{1}{2}
\beta\sum_{\ux,\uy\in\Lambda}J_{\ux,\uy}\sigma_{\ux}\sigma_{\uy}
\right)\times{\rm d}\left(
\otimes_{\ux\in\ZZ^2}\rho_{\rm ss}
\right)(\sigma)
\]
where $\cZ_{\beta,\Lambda}$ is the normalization constant or partition function.
Using the Griffiths inequalities it is easy to construct the infinite volume Gibbs measure $\nu_{\beta}$
at inverse temperature $\beta$ and zero magnetic field as the weak limit of $\nu_{\beta,\Lambda}$ when $\Lambda\nearrow\ZZ^2$.
The corresponding spin correlations $\langle \sigma_{\ux_1}\cdots\sigma_{\ux_n}\rangle=\EE(\sigma_{\ux_1}\cdots\sigma_{\ux_n})$
are well defined and translation invariant.
This choice of boundary condition is called ``free'' in rigorous statistical mechanics and ``Dirichlet'' in constructive quantum field theory.
At the critical temperature, i.e., when $\beta=\beta_{\rm c}=\frac{1}{2}\log(1+\sqrt{2})$, the correlations exhibit power-law decay
at large distances, e.g., $\langle \sigma_{\ux}\sigma_{\uy}\rangle\sim |\ux-\uy|^{-2[\phi]}$
where the symbol $[\phi]$ stands for the {\em scaling dimension}
of the continuum field $\phi$ in the scaling limit. This is also minus the Hurst self-similarity exponent.
For Ising in $d=2$ dimensions, $[\phi]=\frac{1}{8}$. Let $L>1$ be a fixed integer and let $r\in\ZZ$ and for given spin configuration
$\sigma=(\sigma_{\ux})_{\ux\in\ZZ^2}$ define the tempered Schwartz distribution $\Phi_r=L^{r(d-[\phi])}\sum_{\ux\in\ZZ^2}
\sigma_{\ux}\delta_{L^r\ux}$
where $\delta_{L^r\ux}$ is the delta function located at $L^r\ux$. Let $\nu_{\phi,r}$ be the probability measure on $S'(\RR^2)$
obtained as direct image by the previous map from $\nu_{\beta}$. It follows from the combined results of~\cite{CamiaGN,ChelkakHI}
that the scaling limit $\nu_{\phi}$, obtained as weak limit of $\nu_{\phi,r}$ when $r\rightarrow -\infty$,
exists, is non-Gaussian and satisfies (global) CI (in a sense to be explained later).
The articles~\cite{CamiaGN,ChelkakHI}
prove much more and in particular address local CI. However, this program focuses on the ``full plane'' situation
and global CI. Indeed, only the latter is available in dimension $d=3$ and the ultimate
goal of this program is to prove a result analogous to the one mentioned above for the three-dimensional fractional
$\phi^4$ model~\cite{BrydgesMS,AbdesselamCMP}
which conjecturally is the scaling limit of a long-range Ising model~\cite{Mitter,PaulosRRZ}.
The results in~\cite{CamiaGN,ChelkakHI} give rigorous proofs for predictions~\cite{BelavinPZ1,BelavinPZ2}
made by physicists working in {\em conformal field theory} (CFT)
which is a proper subset of {\em quantum field theory} (QFT).
However, from the point of view of theoretical physics, this is rather well understood~\cite{DiFrancescoMS}. 
By contrast, the program outlined in the present article aims at establishing direct
rigorous mathematical contact with very hot physics, i.e., the area of {\em higher-dimensional conformal bootstrap} (HDCB).
The CI of the fractional $\phi_{3}^{4}$ model is the object of a recent prediction~\cite{PaulosRRZ} by researchers in this area.

\section{Long-term view and precise motivating conjectures}

The fractional $\phi_{3}^{4}$ model corresponds to probability measures on distributions $\phi$ in $S'(\RR^3)$ which can be formally written as
\begin{equation}
\frac{1}{\cZ}\exp\left(-\frac{1}{2}\langle \phi,(-\Delta)^{\alpha}\phi\rangle_{L^2}
-\int_{\RR^3}\{g\phi(x)^4+\mu \phi(x)^2\}\ {\rm d}^3x
\right)\ {\rm D}\phi
\label{formal}
\end{equation}
where the fractional power of the Laplacian is $\alpha=\frac{3+\epsilon}{4}$ with $0<\epsilon\ll 1$ and ${\rm D}\phi$
is the nonexistent Lebesgue measure on $S'(\RR^3)$.
The rigorous study of such measures necessitates regularization. Let $\rho_{\rm UV}$ be a mollifier, namely, a smooth function
$\RR^3\rightarrow\RR$ that is invariant by the orthogonal group $O(3)$, has compact support and integral equal to $1$.
Let also $\rho_{IR}$ be a smooth cut-off function, namely, one that is nonnegative, $O(3)$-symmetric,
compactly supported and identically equal to $1$ in a neighborhood of the origin.
Define the non-cut-off covariance $C_{-\infty}$ as the continuous symmetric bilinear form on $S(\RR^3)$ given in Fourier space as
\[
C_{-\infty}(f,g)=\frac{1}{(2\pi)^3}\int_{\RR^3}\frac{\overline{\widehat{f}(\xi)}\widehat{g}(\xi)}{|\xi|^{3-2[\phi]}}\ {\rm d}^3\xi
\]
with $[\phi]=\frac{3-\epsilon}{4}$.
By the Bochner-Minlos Theorem (see, e.g.,~\cite[\S I.2]{Simon2}) there is
a well defined centered Gaussian probability measure $\mu_{C_{-\infty}}$ on $S'(\RR^3)$
with covariance $C_{-\infty}$. For $r\in\ZZ$, define the rescaled mollifier $\rho_{{\rm UV},r}$ by $L^{-3r}\rho_{\rm UV}(L^{-r}\bullet)$.
For $\phi$ sampled according to $\mu_{C_{-\infty}}$, let $\mu_{C_r}$ denote the law of $\phi\ast\rho_{{\rm UV},r}$ and let $C_r$
be the corresponding covariance. The sample paths of $\mu_{C_r}$ are smooth and (with a slight abuse of notation) the kernel
$C_r(x,y)=C_r(x-y)$ of $C_r$ is also smooth. For $s\in\ZZ$ define the rescaled cut-off function $\rho_{{\rm IR},s}$ by
$\rho_{\rm IR}(L^{-s}\bullet)$. Given a {\em bare ansatz} $(g_r,\mu_r)_{r\in\ZZ}$ or rather the germ of such a sequence at $-\infty$, one has
well defined probability measures $\nu_{r,s}$ on $S'(\RR^3)$ given by
\begin{equation}
\frac{1}{\cZ_{r,s}}\exp\left(
-\int_{\RR^3}\rho_{{\rm IR},s}(x)\{g_r:\phi^4:_{C_r}(x)+\mu_r :\phi^2:_{C_r}(x)\}\ {\rm d}^3x
\right)\ {\rm d}\mu_r(\phi)
\label{regularized}
\end{equation}
where Wick ordering $:\phi^2:_{C_r}(x)=\phi(x)^2-C_r(0)$ and
$:\phi^4:_{C_r}(x)=\phi(x)^4-6C_r(0)\phi(x)^2+3C_r(0)^2$ is introduced as a matter of convenience.
The mathematical study of the formal measures (\ref{formal}) amounts to {\em parametrizing the set of all weak limits}
of $\nu_{r,s}$ when
$r\rightarrow -\infty$ and $s\rightarrow \infty$. Since the choice of bare ansatz is infinite-dimensional one may naively expect
that this set is also infinite-dimensional. Conjecturally, because of the {\em universality} phenomenon,
this should instead be a two-dimensional family
(parametrized by the unstable manifold of the trivial RG fixed point as explained in~\cite[p. 6]{Abdesselam2013}).
Of particular interest are isolated points in this set corresponding to
self-similar or scale invariant (see~\cite{Dobrushin1,Dobrushin2}) weak limits
where the resulting random Schwartz distribution $\phi$ satisfies $\lambda^{[\phi]}\phi(\lambda\bullet)
\stackrel{d}{=}\phi(\bullet)$ (doubly!) in the sense of distributions, for suitable exponent $[\phi]$ and all $\lambda>0$.
The high temperature fixed point is Gaussian white noise on $\RR^3$ with $[\phi]=\frac{3}{2}$ and should result when $\mu_r$ is chosen too large
relative to $g_r$. Taking $g_r=\mu_r=0$ gives the previous Gaussian $\mu_{C_{-\infty}}$, i.e., the
fractional Gaussian field $FGF_{\alpha}(\RR^3)$
of~\cite{LodhiaSSW} with $\alpha=\frac{3+\epsilon}{4}$ and $[\phi]=\frac{3-\epsilon}{4}$.
However, there should be a third much more interesting,
i.e., highly non-Gaussian, self-similar measure $\nu_{\phi}$ with the same $[\phi]=\frac{3-\epsilon}{4}$.

\noindent{\bf Conjecture 1:}
{\em Set $[\phi]=\frac{3-\epsilon}{4}$ for $\epsilon$ positive and small enough.
Then there exists a nonempty open interval $I\subset (0,\infty)$ and
a function $\mu_{\rm c}:I\rightarrow\RR$ such that for all $g\in I$, if one sets
$g_r=L^{-r(3-4[\phi])}g$ and $\mu_r=L^{-r(3-2[\phi])}\mu_{\rm c}(g)$, the weak limit $\nu_{\phi}=\lim_{r\rightarrow -\infty}
\lim_{s\rightarrow \infty} \nu_{r,s}$ exists and is non-Gaussian, translation-invariant (stationary), $O(3)$-invariant (isotropic),
scale-invariant (self-similar) with exponent $[\phi]$. Moreover this limit is independent of the choice of $L$, that of $g\in I$, as well as that of the
cut-off functions $\rho_{\rm UV}$ and $\rho_{\rm IR}$.}

Some readers may not be familiar with what ``weak limit'' means in the above conjecture.
Indeed, probability measures in spaces such as $S'(\RR^3)$ are not popular among probabilists who often prefer the use of abstract Wiener spaces
or probability measures on Banach spaces~\cite{Bogachev}. 
The literature on probability theory in spaces of distributions such as $S'(\RR^d)$ or $D'(\RR^d)$ is less well known.
The probabilist's first choice for a reference on the subject often is the pedagogically excellent~\cite{Walsh}. However, it does not cover enough
ground for the needs of constructive QFT and CP. Likewise, white noise theory~\cite{Hida} is too limited in scope since it only treats
one probability measure on $S'(\RR^d)$: white noise!
A gap in the literature is the absence of a very general ``first-quantized''
Kolmogorov-Chentsov (KC) theorem which would roughly say:
if a probability measure $\mu$ on $S'(\RR^d)$ has finite moments of all orders with smooth kernels $S_n(x_1,\ldots,x_n)$
away from the diagonal behaving like
$|S_n(z+\lambda x_1,\ldots,z+\lambda x_n)|\sim\lambda^{-n[\phi]}$ when $\lambda\rightarrow 0$, then $\mu(\cC^{\alpha})=1$ for all
$\alpha<-[\phi]$.
Here the $\cC^{\alpha}$ denote the $\alpha$-H\"older-Zygmund spaces, i.e., Besov spaces $B_{\infty,\infty}^{\alpha}$ featuring in~\cite{BahouriCD,GubinelliIP}, and we are especially interested in the negative $\alpha$ situation.
A first approximation to such an ideal theorem is~\cite[Thm. 1.4.2]{AdlerT} for Gaussian measures.
Much closer to what we have in mind are the versions of such a result given in~\cite[Lemma 9]{MourratW}
as well as~\cite[Thm. 2.7]{ChandraW} and~\cite[Prop. 2.32]{FurlanM}. Note that there is no need for versions and modifications since
we are talking about random distributions~\cite[Ch. III]{Fernique}
rather than the weaker notion of generalized random field~\cite[Ch. II]{Fernique} and, furthermore,
the spaces $\cC^{\alpha}$ are measurable subsets of $S'(\RR^d)$.
The above statement we gave for the ideal KC Theorem is not strictly-speaking correct since one needs a hypothesis on the correlations
{\em as distributions on the diagonal}
which would follow, e.g., if $[\phi]<\frac{d}{2}$, i.e., the kernel is {\em locally integrable}, and correlations are
given by integration against this kernel.
Finally, the ideal KC Theorem should include the case $[\phi]<0$ and the most basic KC Theorem for Brownian motion. This requires the use
of spaces of test functions $S_{k}(\RR^d)$ with vanishing spatial moments up to order $k$ as in~\cite{Dobrushin1}, or non-canonical lifts
for incremental fields as in~\cite{LodhiaSSW}. The simplest example of such a lift is the convention $B_t=0$ at $t=0$ for Brownian motion. Spaces of distributions like $S'(\mathbb{R}^d)$ and $D'(\mathbb{R}^{d})$ are well behaved from the point of view of probability theory. This is largely due to the nuclear property which unfortunately is rather technical.
The author's simplifying point of view is to take as
{\em practical definition} of a nuclear space: a locally convex topological vector space which is isomorphic
to $\mathfrak{s}$, $\mathfrak{s}_0$ or $\mathfrak{s}_0\widehat{\otimes}\mathfrak{s}$. Here $\mathfrak{s}=\mathfrak{s}(\NN)$ is the space of sequences
with faster than power-law decay. One can similarly define $\mathfrak{s}(\NN^d)$ and $\mathfrak{s}(\ZZ^d)$. For instance,
if one used a $\phi^4$-type single site measure ${\rm d}\rho_{\rm ss}(u)\sim\exp(-au^4-bu^2){\rm d}u$, the theory of
superstability (see~\cite[\S 4.7]{Chandra} and references therein) would realize infinite volume lattice measures in $\mathfrak{s}'(\ZZ^2)$
before transfer to $S'(\RR^2)$ for the study of the scaling limit.
As topological vector spaces $\mathfrak{s}\simeq \mathfrak{s}(\NN^d)\simeq \mathfrak{s}(\ZZ^d)\simeq S(\RR^d) $~\cite{Simon1}.
The space $\mathfrak{s}_0=\oplus_{\NN}\RR$ is that of almost finite sequences with the finest locally convex topology, i.e., the one generated
by the uncountable set of all possible seminorms. The Schwartz-Bruhat spaces~\cite{Bruhat,Osborne} $S(\QQ_p^d)$ are isomorphic to
$\mathfrak{s}_0$. Finally $\mathfrak{s}_0\widehat{\otimes}\mathfrak{s}=\oplus_{\NN}\mathfrak{s}$ with the topology defined by the set of all seminorms
that are continuous on individual summands. The adelic Schwartz-Bruhat space $S(\mathbb{A}_{\QQ})$
is isomorphic to $\mathfrak{s}_0\otimes\mathfrak{s}$ but so is the {\em classical
space of test functions} $D(U)$, for $U$ a nonempty open set
in $\RR^d$. Explicit Schauder bases realizing the last Valdivia-Vogt isomorphism~\cite{Valdivia,Vogt} were not known until
recently~\cite{Bargetz}.
In the statement of Conj. 1, $S'(\RR^d)$ is to be understood as a {\em topological space} with the strong topology. The $\sigma$-algebra is the Borel $\sigma$-algebra
for this topology.
One can also use the weak-$\ast$ topology, but the Borel $\sigma$-algebra would be the same and it would also 
equal the cylinder $\sigma$-algebra (that makes the maps $\phi\mapsto \phi(f)$ measurable)~\cite{Fernique,Becnel}.
More importantly, {\em weak convergence of probability measures is defined as usual via the convergence of expectations of bounded continuous functions}.
By a remarkable theorem of Fernique~\cite[Thm. III.6.5]{Fernique}, this definition is insensitive to the choice of topology between
the weak-$\ast$ and strong ones
in the definition of bounded {\em continuous} functions.
While $S'(\RR^d)$ is non-metrizable and thus not Polish (the popular tools in~\cite{Billingsley} fail), probability theory on this
space is in some sense {\em nicer} than in Banach spaces. Prokhorov's Theorem, Bochner's Theorem and the L\'evy Continuity Theorem
all apply to $S'(\RR^d)$ and their statements are {\em identical} to their finite-dimensional analogues (see~\cite{Fernique,Walsh}).

A probability measure $\mu$ on $S'(\RR^d)$ has the {\em moments of all orders} (MAO) property if $\phi\mapsto\phi(f)$ is
in $L^p(S'(\RR^d),\mu)$ for all $f\in S(\RR^d)$ and all $p\in [1,\infty)$. By~\cite[Cor. II.2.5]{Fernique}, the moments
$S_n(f_1,\ldots,f_n)=\langle \phi(f_1)\cdots\phi(f_n) \rangle=\int \phi(f_1)\cdots\phi(f_n)\ {\rm d}\mu(\phi)$
are automatically continuous as $n$-linear forms on $S(\RR^d)$.
By the Schwartz-Grothendieck Nuclear Theorem (see~\cite[Thm. 6]{Simon1} for an easy proof), the moments $S_n$ can
be seen as distributions, i.e., elements of $S'(\RR^{nd})$. 
An MAO measure $\mu$ is {\em determined by correlations} (DC) if the only MAO measure with the same sequence of moments $S_{n}$ as $\mu$
is $\mu$ itself.
A DC measure will be called {\em determined by pointwise correlations} (DPC) if: 1) for all $n$, the $S_n\in S'(\RR^{nd})$ have singular support
(see, e.g.,~\cite{BrouderDH}) contained in the
big diagonal ${\rm Diag}_n=\{(x_1,\ldots,x_n)\in\RR^{nd} | \exists i\neq j, x_i=x_j\}$ (this property uniquely defines the
smooth pointwise functions
$S_n(x_1,\ldots,x_n)=\langle \phi(x_1)\cdots\phi(x_n)\rangle$ on $\RR^{nd}\backslash {\rm Diag}_n$); 2) these pointwise correlations are
$L^{1,{\rm loc}}$ on ${\rm Diag}_n$; and 3) for all $n$ and all test functions $f_1,\ldots,f_n\in S(\RR^d)$ one has
\[
\langle \phi(f_1)\cdots\phi(f_n)\rangle=\int_{\RR^{nd}\backslash {\rm Diag}_n}
\langle \phi(x_1)\cdots\phi(x_n)\rangle f(x_1)\cdots f(x_n)\ {\rm d}^dx_1\cdots
{\rm d}^dx_n
\]
where on the left one has an honest expectation with respect to $\mu$ whereas on the right, the expectation-like
notation $\langle \phi(x_1)\cdots\phi(x_n)\rangle$ is merely {\em symbolic}.

\noindent{\bf Conjecture 2:} {\em The measure $\nu_{\phi}$ from Conj. 1 is DPC.}

Note that pointwise correlations must satisfy $\langle \phi(\lambda x_1)\cdots\phi(\lambda x_n)\rangle
=\lambda^{-n[\phi]}\langle \phi(x_1)\cdots\phi(x_n)\rangle$ as a result of self-similarity.
If one keeps in mind translation invariance, then local integrability, e.g., for $n=2$ amounts to $d-2[\phi]>0$ which is luckily satified
for $d=3$ and $[\phi]=\frac{3-\epsilon}{4}$ with $\epsilon$ small.

\noindent{\bf Conjecture 3:} {\em The measure $\nu_{\phi}$ is Osterwalder-Schrader (OS) positive (see~\cite[\S 6.1]{GlimmJ}).}

Thus the measure should satisfy all the OS axioms and should provide a nontrivial Euclidean QFT in the traditional sense of constructive QFT.
Although the last conjecture does not pertain to probability per se, it has important consequences with relation to probability.
For instance it would imply {\em real analyticity} for the pointwise correlations $\langle \phi(x_1)\cdots\phi(x_n)\rangle$
(see~\cite[Cor. 19.5.6]{GlimmJ}). OS positivity can also help prove bounds on correlations (see~\cite{FrohlichSS} for a famous example
and~\cite{Shlosman} for a review).
Since the recent work of Hairer~\cite{Hairer,Hairer2} (see also~\cite{CatellierC,Kupiainen})
on the stochastic quantization SPDE for the $\phi_3^4$ model, there has been renewed attention devoted to the latter.
However, in~\cite{Jaffe}, the author expresses some criticism for this approach seen as a way to construct the $\phi_3^4$
model and show that it satisfies the OS axioms.
Such criticism is based on the premise that it would essentially be impossible to prove OS positivity for a nontrivial limit
(such as $\nu_{\phi}$) if it did not hold for the approximants (such as $\nu_{r,s}$). This premise may be too pessimistic.
Indeed, the author noticed a few years ago (talk at the conference
``Rigorous Quantum Field Theory in The LHC Era'', Vienna, 2011) that all one needs is
a {\em partial} OS positivity for approximants that becomes full OS positivity
in the limit. Namely, one only needs positivity for observables that are ``$L^r$ away'' from the reflection hyperplane. This is what
regularization on a lattice of mesh $L^r$ does!
Thus a strategy for proving Conj. 3 could be as follows. Let $\theta (x_1,\ldots,x_d)=\theta (-x_1,x_2,\ldots,x_d)$, $\theta f(x)=f(\theta x)$,
$\theta\phi(f)=\phi(\theta f)$ and $\theta F(\phi)=F(\theta\phi)$
denote the reflection respectively for points, test functions, distributional
fields and measurable functionals.
OS positivity would follow if one can show $\int \overline{\theta F(\phi)} F(\phi){\rm d}\nu_{\phi}(\phi)\ge 0$
for all $F$'s given by a polynomial in the ``coordinates'' $\phi(f_1),\ldots,\phi(f_n)$ with test functions $f_1,\ldots,f_n$
having compact support in the open
half-space $(0,\infty)\times\RR^{d-1}$. Once these functions are fixed then there is a
lower bound $a>0$ for the distance between these supports and
the reflection hyperplane. Now replace the couplings $g_r$ and $\mu_r$ in (\ref{regularized}) by space-dependent functions
$\widetilde{g}_r(x)$ and $\widetilde{\mu}_r(x)$ which vanish in a corridor of thickness $\sim L^r$ around the reflection hyperplane, are
essentially equal to the previous uniform values $g_r$, $\mu_r$ further away,
and may eventually be bigger at the edges of the corridor in order to compensate
for the inside loss. Then the analogue of $\int \overline{\theta F(\phi)} F(\phi){\rm d}\nu_{r,s}(\phi)$ can be written
as $\int \overline{\theta F(\phi\ast\rho_{{\rm UV},r}) e^{-\theta\widetilde{V}_r(\phi\ast\rho_{{\rm UV},r})}}
F(\phi\ast\rho_{{\rm UV},r})e^{-\widetilde{V}_r(\phi\ast\rho_{{\rm UV},r})}
{\rm d}\mu_{C_{-\infty}}(\phi)$
for a suitable functional $\widetilde{V}_r$. Roughly speaking, the latter lives in $(L^r,\infty)\times\RR^{d-1}$.
The needed positivity then is a result of OS positivity for the {\em non-cut-off} Gaussian $\mu_{C_{-\infty}}$,
the size $\sim L^r$ for the compact support of $\rho_{{\rm UV},r}$ and taking $L^r\ll a$. Hence,
the problem is reduced to showing that this evanescent corridor modification has no effect on the weak limit $\nu_{\phi}$.
One thus needs a {\em robustness} result with respect to a perturbation of the {\em interaction}
(the part inside the exponential in
(\ref{regularized})) corresponding to a hypersurface {\em defect} in the middle of the bulk.
Note that a similar robustness is needed in Conj. 1 in order to deal with the smooth fall-off of $\rho_{{\rm IR},s}$.
While $O(d)$ invariance in Conj. 1 is trivial as stated, this is not the case for translation invariance.
One needs similar robustness with respect to {\em boundary} effects due to the symmetric difference of two large spheres.
The latter should be negligible compared to the bulk since the shape satisfies the Van Hove criterion
(see~\cite[Def. 2.1.1]{Ruelle}) for thermodynamic limits, but one has to prove it.

Given two nonempty simply connected open sets $U$ and $V$ in $\RR^d$, a $C^k$-diffeomorphism $f:U\rightarrow V$,
with $k\ge 1$ finite or infinite,
is called {\em conformal} if there exists a function $\tau:U\rightarrow \RR$ such that $\forall x\in U$,
$e^{-\tau(x)} {\rm D}_x f\in O(d)$. The condition is void for $d=1$. If $d=2$, then $f(z)$ must be of the form $g(z)$ or $g(\overline{z})$ for
some holomorphic map $g$.
If one does not fix the domain $U$ and target $V$ this gives a rich infinite-dimensional flexibility
(as in the titles of~\cite{BelavinPZ1,BelavinPZ2});
but if one does
then one is reduced to $SL_2(\CC)$ or $SL_2(\RR)$ M\"{o}bius symmetry.
Finally, if $d\ge 3$ then by a theorem of Liouville (see, e.g.,~\cite{Hartman}) $f$ must be the restriction of an element of the M\"{o}bius
group $\cM(\RR^d)$ of global conformal maps.
This group can be defined as the group of bijective transformations of the ``sphere'' $\widehat{\RR^d}=\RR^d\cup\{\infty\}$
generated by Euclidean isometries, scaling transformations $x\mapsto \lambda x$ with $\lambda>0$ and the unit-sphere inversion
$J(x)=|x|^{-2} x$. By definition, the first two preserve the point at infinity while $J$ exchanges the origin and $\infty$.
Another more elegant definition of $\cM(\RR^d)$ can be given using the {\em absolute cross-ratio}. In the convention of~\cite{Ahlfors}
(there is no standard choice in the literature) it is given by ${\rm CR}(x_1,x_2,x_3,x_4)=(|x_1-x_3|\ |x_2-x_4|)/(|x_1-x_4|\ |x_2-x_3|)$
for quadruples in $(\RR^d)^4\backslash{\rm Diag}_4$.
It can be extended by continuity ($\widehat{\RR^d}$ is seen as the one-point compactification of $\RR^d$) to
$\widehat{\RR^d}\backslash \widehat{{\rm Diag}}_4$ where $\widehat{{\rm Diag}}_4$ is the new big diagonal, by omitting factors which
contain $\infty$.
Then $\cM(\RR^d)$ {\em is} the group of bijective transformations which leave the cross-ratio invariant
(see~\cite[Thm. 4]{Wilker} or~\cite[Thm. 4.3.1]{Ratcliffe}). 

\noindent{\bf Conjecture 4:}
{\em The measure $\nu_{\phi}$ satisfies global CI. Namely, the pointwise correlations for $\nu_{\phi}$
obey
\begin{equation}
\langle \phi(x_1)\cdots\phi(x_n)\rangle=
\left(
\prod_{i=1}^{n} |J_f(x_i)|^{\frac{[\phi]}{3}}
\right)\times
\langle \phi(f(x_1))\cdots\phi(f(x_n))\rangle
\label{pointwiseCI}
\end{equation}
for all $f\in\cM(\RR^3)$, and all collections of distinct points in $\RR^3\backslash\{f^{-1}(\infty)\}$.
Here, $J_f(x)={\rm det}({\rm D}_x f) $ denotes the Jacobian determinant of $f$ at $x$.}

This is the analogue of the self-similarity propery stated after Conj. 2 where now
the rescaling factor $\lambda$ is allowed to be space-dependent.
Conj. 4 is supported by the arguments in~\cite{PaulosRRZ}.
CI for the fractional $\phi_{3}^{4}$ model thus joins similar physically motivated conjectures in 3D for
the self-avoiding walk~\cite{Kennedy},
percolation~\cite{GoriT} or the short-range Ising model (see, e.g.,~\cite{DelamotteTW}). 
The 2D Ising analogue of Conj. 4 immediately follows from the explicit formula
proved in~\cite[Eq. 1.6]{ChelkakHI} for the full plane case
(see also~\cite[Thm. 4]{Dubedat}). Note that the argument in~\cite{PaulosRRZ} is quite subtle since the fractional
$\phi_{3}^{4}$ model is {\em nonlocal}
and one is missing the most important ingredient of textbook CFT~\cite{DiFrancescoMS}, namely, a local energy-momentum tensor.
Instead,~\cite{PaulosRRZ} makes fundamental use of the point of view known in the physics literature as
${\rm AdS}_{d+1}/{\rm CFT}_d$ correspondence and holography (see, e.g.,~\cite{Witten,Kaplan})
which involves a notion of extension from the boundary to the bulk as in~\cite{CaffarelliS}.
Indeed, $\widehat{\RR^d}$ can be identified with the unit sphere $\mathbb{S}^d$ inside $\RR^{d+1}$. This sphere is the topological
boundary of the open unit ball $\mathbb{B}^{d+1}$. Although not a conformal invariant, the metric induced from the Euclidean one on
$\RR^{d+1}$ is a reasonable one to put
on $\mathbb{S}^d$. As for $\mathbb{B}^{d+1}$, the good metric to use is the hyperbolic one given by ${\rm d}s=\frac{2|{\rm d}x|}{1-|x|^2}$.
The AdS/CFT point of view is based on the classical theorem in geometry which establishes a one-to-one correspondence between
elements of $\cM(\RR^d)$ acting on the boundary $\widehat{\RR^d}\simeq \mathbb{S}^d$ and hyperbolic {\em isometries}
acting on $\mathbb{B}^{d+1}$
(see, e.g.,~\cite[Thm. 4.5.2]{Ratcliffe}). Instead of using the previous {\em conformal ball model} of hyperbolic geometry
one can also use the {\em upper half-space model} where $\widehat{\RR^d}$ is identified with the hyperplane $\RR^d\times\{0\}\subset\RR^{d+1}$
together with the point at infinity ``$(0,\ldots,0,\infty)$''. The hyperplane is equipped with the usual Euclidean metric and appears
as the boundary (in $\RR^{d+1}$) of the upper half-space $\mathbb{H}^{d+1}=\RR^d\times (0,\infty)$.
The latter is equipped with the hyperbolic metric
${\rm d}s=|{\rm d}x|/x_{d+1}$. The spaces $\mathbb{B}^{d+1}$ or $\mathbb{H}^{d+1}$ are what physicists
call the {\em Euclidean Anti-de Sitter space} ${\rm AdS}_{d+1}$
(the word ``Euclidean'' is {\em not} in the sense of geometry but in the sense of the Euclidean vs. Lorentzian/Minkowskian distinction
in axiomatic and constructive QFT).

The article~\cite{PaulosRRZ} not only covers the 3D fractional $\phi^4$ model but, more generally, critical scaling limits
of long-range Ising models in any dimension.
Indeed, instead of using the mollifier $\rho_{\rm UV}$ one can use a lattice regularization, as in~\cite{BrydgesFS}, in order to make
sense of (\ref{formal}). The fractional Laplacian $(-\Delta)^{\alpha}$ with $\alpha=\frac{3+\epsilon}{4}$
would be replaced by the lattice analogue
$(-\Delta_{L^r\ZZ^3})^{\alpha}$ adapted to the fine lattice with mesh $L^r$. This is a rescaled version of the unit lattice fractional
Laplacian $(-\Delta_{\ZZ^3})^{\alpha}$ with off-diagonal matrix elements $-J_{\ux,\uy}$. It is an easy exercise to check that $J_{\ux,\uy}>0$
(the model is {\em ferromagnetic})
and $J_{\ux,\uy}\sim |\ux-\uy|^{-(d+\sigma)}$ at long distance, with $d=3$ and
$\sigma=\frac{3+\epsilon}{2}$. Long-range models were rigorously
studied in $d=1$ by Dyson~\cite{Dyson1,Dyson2} and for general $d$ and $\sigma$ in~\cite{AizenmanF} where the existence of a critical
inverse temperature $\beta_{\rm c}$ is shown. The corresponding critical scaling limits were the object of non mathematically rigorous investigations
by physicists~\cite{FisherMN,Sak1,Sak2}. For $d=3$, as one increases $\sigma$ one should see three regimes: mean-field, intermediate, and
short-range (see~\cite[Fig. 1]{PaulosRRZ}).
The $\epsilon\sim 1$ regime of the BMS model where one would transition into the more mysterious scaling limit for the 3D short-range Ising
model has been the object of recent controversy among physicists~\cite{Picco,BlanchardPR,AngeliniPR,BrezinPR,BehanRRZ}.
This program pertains to the other more mathematically tractable
end $0<\epsilon\ll 1$ of the intermediate region, where one is close to the mean-field regime yet below the upper critical dimension.
Note that if one uses the lattice regularization for (\ref{formal}) and makes the same choices as in Conj. 1 for the bare ansatz $g_r,\mu_r$ 
then one is taking the scaling limit of a {\em fixed} critical lattice measure for the long-range model (similarly to the discussion
of~\cite{CamiaGN,ChelkakHI}
in relation to the 2D short-range Ising model).

\noindent{\bf Conjecture 5:}
{\em The measure $\nu_{\phi}$ can be obtained as a critical lattice scaling limit starting with 1) a single site measure $\rho_{\rm ss}$
of $\phi^4$ type as well as with 2) $\rho_{\rm ss}=\frac{1}{2}(\delta_{-1}+\delta_{1})$.
}

We have nothing to say at this point about Conj. 5. Indeed Part 1) of the latter amounts to robustness with respect to the choice of
regularization or perturbations of the Gaussian measure (independence from the choice of $\rho_{\rm UV}$ in Conj. 1).
This is very hard to do with rigorous RG methods (see however~\cite{BauerschmidtBS4} where a particular instance of this
issue was overcome!). {\em The leitmotiv of this program is robustness with respect to
changes in the interaction potential} (e.g., the choice of $\rho_{\rm IR}$ in Conj. 1).
Part 2) of Conj. 5 seems even more difficult.
This conjecture in the physics literature
~\cite[\S 7.4.7]{DiFrancescoMS} goes under the heading of effective Ginzburg-Landau (or multi-critical model)
description of the unitary minimal models $\cM(m+1,m)$ of CFT~\cite{Zamolodchikov}.
If one uses block-spins to implement the RG, then one would sum unit-lattice spins $\sigma_{\ux}$ over blocks of size $L$ and multiply
by $L^{-(d-[\phi])}$. The usual coarse-graining for the central limit theorem is when $[\phi]=\frac{d}{2}$ as in~\cite{Newman1,Newman2}
or~\cite{BreuerM}. If one repeats the operation $N$ times then the spacing between spin values is $\sim L^{-N(d-[\phi])}\rightarrow 0$
whereas the extreme values go as $L^{N[\phi]}\rightarrow\infty$ so the spin distribution looks like a law on $\RR$ with continuous density.
If one starts with $\rho_{\rm ss}=\frac{1}{2}(\delta_{-1}+\delta_{1})$ it is therefore intuitively reasonable to expect the same scaling limit
as if one started with a double-well $\phi^4$-type $\rho_{\rm ss}$ (see also~\cite[\S 12.3]{Kupiainen2014}).
Finally, with regard to the $\phi_2^4$/Ising conjectural equivalence in 2D, note that massless measures $\nu_{m=0}$ are known to
exist (also for $\phi_3^4$) from the works~\cite{McBryanR,GlimmJ1977} (see also~\cite[\S9]{BrydgesFS}).
However, these are {\em not} self-similar and thus $\nu_{m=0}\neq\nu_{\phi}$ since the bare coupling $g_r$ is kept fixed with respect to
the UV cut-off $r$ (see~\cite[\S4]{BrydgesFS}). One would need to take the large distance scaling limit (see~\cite[Eq. 3.1]{Dobrushin2})
of $\nu_{m=0}$
in order to recover $\nu_{\phi}$. From the RG perspective, $\nu_{m=0}$ corresponds to a two-sided RG trajectory
starting at a Gaussian fixed point and ending at a nontrivial fixed point, just like the one constructed in~\cite{AbdesselamCMP}
for the 3D fractional $\phi^4$ model.

The article~\cite{PaulosRRZ} belongs to HDCB which has known tremendous activity in the last few years.
For a long time, it was believed that conformal symmetry in $d\ge 3$ was too poor to produce precise predictions for critical
exponents as 2D CFT did; but the situation changed!
The, soon to explained, {\em operator product expansion} (OPE) is the main conceptual foundation for
HDCB~\cite{Rychkov,ElShowketal} as well as for 2D CFT~\cite{BelavinPZ1,BelavinPZ2}.
The OPE is expected to hold in any QFT conformal or not, and it gives precise predictions for the
asymptotic behavior of pointwise correlations near ${\rm Diag}_n$.
Supplemented with CI and the associative property (crossing symmetry), the OPE essentially allows one to
compute the correlations
of the CFT at hand. 
Recent progress was made possible
by better formulas for conformal blocks~\cite{DolanO1,DolanO2}, and the discovery~\cite{Rattazzietal}
of a universal bound $[\phi^2]\le f([\phi])$ relating the scaling dimension of the field $\phi$, with law $\nu_{\phi}$,
to the scaling dimension $[\phi^2]$ of a suitably renormalized pointwise square $N[\phi^2]$. The 2D Ising model was 
found to be a ``kink'' which saturates this bound at the point $([\phi],[\phi^2])=(\frac{1}{8},1)$.
Similarly to thinking of $\phi$ as the scaling limit of the spin field, a good way to think of $N[\phi^2]$
is as the scaling limit of the energy field whose correlations were studied in~\cite{Hongler,HonglerS}.
It was observed (see~\cite{Rychkov}) that a similar positioning at a kink happens for 3D Ising 
and this could give a clue to the pair $([\phi],[\phi^2])$ for 3D Ising.
Some of the most precise theoretical predictions for 3D Ising critical exponents are in~\cite{SimmonsD}, namely: $[\phi]\simeq
0.518151$ and $[\phi^2]\simeq 1.41264$.
The author previously thought that such HDCB
predictions were the result of an appealing but unjustified ``kink hypothesis'' for 3D Ising, but in fact
they rest on much more solid arguments using OS positivity, the OPE and crossing symmetry applied to the {\em mixed} four-point functions
of the fields $\phi$ {\em and} $N[\phi^2]$~\cite{KosPS}.
For the fractional $\phi_{3}^{4}$ model in $\RR^3$ one expects $[\phi^2]-2[\phi]=\frac{\epsilon}{3}+o(\epsilon)$
(see, e.g.,~\cite[Eq. 2.18]{PaulosRRZ}). This is {\em exactly}
the relation proved in~\cite{AbdesselamCG1} for the hierarchical fractional $\phi_{3}^{4}$ model.
Such composite field {\em anomalous scaling} is reminiscent of {\em multifractality} (see, e.g.,~\cite{AllezRV}).
For an ordinary random field $\phi(x)$, with moments satifying
$\EE |\phi(x)-\phi(y)|^m\sim |x-y|^{-[\phi^m]}$, multifractal behavior can be described as a nonlinear variation of the scaling exponent
$[\phi^m]$ in terms of the power $m$. When the Hurst exponent $-[\phi]$ is positive, the process
is not a generalized one, and the power is a ``true power'',
then such behavior is ruled out by strict self-similarity. For the fractional $\phi_{3}^{4}$ model the strictly self-similar field $\phi$ sampled according
to $\nu_{\phi}$
is a generalized one and the square
$N[\phi^2]$ is {\em not} a naive
product of the field with itself. In addition to additive renormalizations as in the usual Wick product construction~\cite{DaPratoT} or
the recent work~\cite{Hairer}, one needs a nontrivial {\em multiplicative}
renormalization which is responsible for this ``multifractal''
relation $[\phi^2]>2[\phi]$.
It is remarkable for the {\em simplified} model in~\cite{AbdesselamCG1} to produce, if one boldly keeps the $O(\epsilon)$ terms and extrapolates
to $\epsilon\sim 1$ (3D short-range Ising), a value $1/3=0.33\ldots$ for the scaling anomaly $[\phi^2]-2[\phi]$ that is quite close
to the value $0.3763\ldots$ deduced from~\cite{SimmonsD}.

\noindent{\bf Conjecture 6:}
{\em There exists $[\phi^2]$, a function of $\epsilon$ only, satisfying
$[\phi^2]-2[\phi]=\frac{\epsilon}{3}+o(\epsilon)$ where $[\phi]=\frac{3-\epsilon}{4}$ and $\epsilon$ is supposed to
be small for which the following is true.
Let $\phi_{r,s}$ be a random (smooth) distribution in $S'(\RR^3)$ sampled according to $\nu_{r,s}$ defined by (\ref{regularized}).
Define a new smooth field $N_{r,s}[\phi^2]$ which is given deterministically in terms of $\phi_{r,s}$ by
$N_{r,s}[\phi^2](x)=L^{-r([\phi^2]-2[\phi])}\left[\phi_{r,s}(x)^2-\EE (\phi_{r,s}(x)^2)\right]$
where the expectation is with respect to ${\rm d}\nu_{r,s}(\phi_{r,s})$.
Then the pair of random variables $(\phi_{r,s},N_{r,s}[\phi^2])\in S'(\RR^3)^2$ converges in distribution to a pair $(\phi,N[\phi^2])$
with joint law $\nu_{\phi\times\phi^2}$
such that ${\rm var}(N[\phi^2](f))>0$ for some test function $f$.
}

Clearly, the non-Gaussian measure $\nu_{\phi}$ must be the first marginal of $\nu_{\phi\times\phi^2}$.

\noindent{\bf Conjecture 7:}
{\em The joint law $\nu_{\phi\times\phi^2}$ is DPC as well as translation, rotation, scale and conformal
invariant as in the previous conjectures 1, 2 and 4.}

Of course, $[\phi]$ needs to be replaced by $[\phi^2]$ where appropriate in the statements obtained by trivially adapting those of the previous
conjectures. 
Also, the symmetries act diagonally: the same geometric transformation is applied to $\phi$ and $N[\phi^2]$.
A result of Conj. 7 is to give access to all mixed pointwise correlations
$\langle\ \phi(x_1)\cdots\phi(x_n)\ N[\phi^2](y_1)\cdots N[\phi^2](y_m)\ \rangle$
for the measure $\nu_{\phi\times\phi^2}$ on $S'(\RR^3)^2$. Note that for 2D Ising, the 1st marginal
$\nu_{\phi}$ has been obtained in~\cite{CamiaGN}
but the existence of the joint law $\nu_{\phi\times\phi^2}$ or even just the 2nd marginal $\nu_{\phi^2}$ is more problematic
because $[\phi^2]=1$ and the covariance kernel fails to be $L^{1,{\rm loc}}$. Nevertheless, the mixed
pointwise correlations still make sense and they
should be the scaling limits of mixed spin-energy correlations.
For the fractional $\phi_{3}^{4}$  the situation is better because of the relation $[\phi^2]-2[\phi]=\frac{\epsilon}{3}+o(\epsilon)$
which implies the $L^{1,{\rm loc}}$ property $3-2[\phi^2]>0$ when $\epsilon$ is small.

\noindent{\bf Conjecture 8:}
{\em
There exists positive constants $c_{1,1;0}$ and $c_{1,1;2}$ such that
$\phi(x)\phi(y)=c_{1,1;0} |x-y|^{-2[\phi]}+c_{1,1;2}|x-y|^{[\phi^2]-2[\phi]}
N[\phi^2](x)+o(|x-y|^{[\phi^2]-2[\phi]})$ when $y\rightarrow x$,
inside arbitrary mixed pointwise correlations.
}

The meaning of `inside correlations' should be clear. It implies, e.g.,
$\langle\phi(x)\phi(y)\phi(z_1)\cdots\phi(z_n)\rangle=
c_{1,1;0} |x-y|^{-2[\phi]}\ \langle\phi(z_1)\cdots\phi(z_n)\rangle
+c_{1,1;2}|x-y|^{[\phi^2]-2[\phi]}\langle N[\phi^2](x)\phi(z_1)\cdots\phi(z_n)\rangle
+o(|x-y|^{[\phi^2]-2[\phi]})$
when $y\rightarrow x$ for fixed $(x,z_1,\ldots,z_n)\in\RR^{(n+1)d}\backslash{\rm Diag}_{n+1}$.
The simplest example of the OPE {\em is} Conj. 8.
Note that complete dependence of pairs of random variables does not necessarily survive convergence in distibution. Stronger notions than
weak limits of probability measures would be needed in order to preserve that property (see, e.g.,~\cite{SiburgS}). While 
$N_{r,s}[\phi^2]$ is a function of $\phi_{r,s}$ it is not clear from Conj. 6 that the ``energy'' field
$N[\phi^2]$ should be a deterministic function, i.e., a pathwise renormalized square of the ``spin'' field $\phi$.

\noindent{\bf Conjecture 9:} {\em Let
$\widetilde{N}_r[\phi^2](x)=L^{-r([\phi^2]-2[\phi])}\left[
(\phi\ast\rho_{{\rm UV},r})(x)^2-\EE \left((\phi\ast\rho_{{\rm UV},r})(x)^2\right)
\right]$
where $\phi$ is sampled according to $\nu_{\phi}$ and
the expectation is with respect to that measure.
Then for all test functions $f$, 
the smeared quantities
$\int_{\RR^3} \widetilde{N}_r[\phi^2](x)\ f(x)\ {\rm d}^3x$, built from the smooth field $\widetilde{N}_r[\phi^2](x)$,
converge in every $L^p(S'(\RR^3),\nu_{\phi})$, for $p\in [1,\infty)$, and almost surely to a limit $\widetilde{N}[\phi^2](f)$.
Moreover, there exists a Borel measurable map ${\rm Sq}:S'(\RR^3)\rightarrow
S'(\RR^3)$
such that 1) for all $f$, ${\rm Sq}(\phi)(f)=\widetilde{N}[\phi^2](f)$, $\nu_{\phi}$-a. s. in $\phi$; and 2)
the law of $(\phi,\lambda{\rm Sq}(\phi))$ is $\nu_{\phi\times\phi^2}$, for a suitable constant $\lambda>0$.
}

As is well known, one cannot define the pointwise product of distributions in general~\cite{Schwartz}.
One could argue that the whole history of stochastic calculus
can be seen as an effort to circumvent this no-go theorem, albeit in a stochastic sense or, even better, in an almost sure or pathwise sense.
This was nicely recalled in~\cite[pp. 5--6]{GubinelliIP} starting with the simplest case: It\^{o}'s integral.
Indeed, when considering $\int f(B_t) {\rm d}B_t$ with $f$ smooth, or rather $\int f(B_t)\dot{B}_t {\rm d}t$ the issue
is not so much integration
as it is pointwise multiplication of $f(B_t)$ (of regularity $\cC^{\alpha}$, $\alpha<\frac{1}{2}$) with $\dot{B}_t$
(of regularity $\cC^{\beta}$, $\beta<-\frac{1}{2}$).
Deterministic multiplication theorems, for instance based on Bony paraproducts~\cite{BahouriCD},
require $\alpha+\beta>0$ which one barely misses in the case of Brownian motion. It\^{o}'s theory circumvents this issue
by using a stochastic construction. Later (see, e.g.,~\cite{Follmer}) one could do this construction in a pathwise or almost
sure sense. Typical pathwise approaches amount to identifying a set $\cG$ of good paths for which one can do the wanted construction
deterministically. Then, one brings in probability by showing that sample paths are a. s. in $\cG$.
Conj. 9 gives a similar pathwise construction of the pointwise square of $\phi$, except the set $\cG$ is implicitly defined and
the separation with probability is not as clear-cut.
Still, Conj. 9 is clearly related to pathwise approaches for products of random distributions as in~\cite{Hairer,GubinelliIP}
which took their inspiration from the theory of rough paths~\cite{Lyons}.
Going from Conj. 8 to Conj. 9 can be done using what we would like to call a ``second-quantized'' KC Theorem as in the recent article~\cite{AbdesselamNew2}.

\noindent{\bf Conjecture 10:}
{\em The field $\phi$ sampled according to $\nu_{\phi}$ is not subordinated to a Gaussian field in the sense
of~\cite{Dobrushin1,Major}
if one only allows finite order Wiener chaos representations.}

Conj. 10 is the explanation of the word ``highly'' used when saying that $\nu_{\phi}$ is a highly non-Gaussian measure.
In constructive QFT one has different grades of nontriviality which can be ranked by
increasing strength as: NT1) the Euclidean probability measure is non-Gaussian; NT2) the corresponding
quantum field obtained by OS reconstruction is not in the Borchers class (see, e.g.,~\cite[\S 4.6]{StreaterW}) of a generalized free field;
NT3) the $S$-matrix is not the identity; NT4) the $S$-matrix exhibits particle production/destruction.
Only NT4 gives Einstein's $E=mc^2$ its full punch!
Note that a Wick square of a free field would be non-Gaussian, i.e., one can have NT1 but not NT2 (see~\cite{Read,Rehren}) and thus
a trivial $S$-matrix~\cite[Cor. 12.2]{BogolubovLOT}.
Conj. 10 is related to property NT2 and it follows easily from Conj. 9 and the existence of an anomalous dimension $[\phi^2]-2[\phi]>0$
from Conj. 6.
Indeed, Conj. 9 says $N[\phi^2]$ is a functional of $\phi$. If the latter was itself a functional,
say $:\psi^k:$, of a Gaussian field $\psi$,
then by the easy analogues of Conj. 8 and 9 for $\psi$ one should be able to rule out a relation such as
$[\phi^2]-2[\phi]=\frac{\epsilon}{3}+o(\epsilon)$ and get a contradiction.
Interestingly, Conj. 10
leaves the door open to an infinite chaos-order representation. In some sense, the Coulomb gas formalism of CFT~\cite[Ch. 9]{DiFrancescoMS}
supports that idea: correlations of nontrivial CFTs like 2D Ising can be obtained using infinite chaos-order fields $:e^{i\alpha \psi(x)}:$
where $\psi$ is the massless free field, modulo suitable screening charges. This ties in with the multifractal story told earlier.
The first precise mathematical model for such behavior, following earlier insights by Kolmogorov and Yaglom, is Mandelbrot's cascade
(see~\cite{AllezRV} and references therein)
initially defined by an RG-like fixed point equation but which can also be seen as a random measure with density (with respect to
Lebesgue measure) given by an infinite chaos-order field $:e^{\alpha \psi(x)}:$ with $\psi$ similar to a hierarchical $\QQ_{2}^{1}$ field
(with $[\phi]=0$, i.e., log-correlated as in~\cite{DuplantierRSV}) in the notations of the next section.
Note that for the fractional $\phi_{3}^{4}$ it is not clear how one can even formulate NT3 and NT4, so Conj. 10 might be the best one can do in the near future.
Indeed, there is no reasonable notion of particles even for the free QFT corresponding to $\mu_{C_{-\infty}}$
studied, e.g., in~\cite{DutschR1,DutschR2}.
The current frontier regarding scattering theory for massless QFTs does not concern relativistic ones (as the one obtained from
$\nu_{\phi}$ by OS reconstruction) but non-relativistic models (see~\cite{ChenFP,Dybalski,DeRoeckK} for recent work in the area).

From the probabilistic point of view, Conj. 4 using moments and their pointwise limits
may seem as a poor way to formulate CI for a probability measure such as $\nu_{\phi}$.
A more satisfactory approach would be to define $T_f\phi$ resulting from the covariant action of a conformal map $f\in\cM(\RR^d)$
on a distribution $\phi\in S'(\RR^d)$. As far as we know, this issue was not addressed in the 2D CP literature, perhaps because
of the emphasis on local rather than global CI.
If $f:U\rightarrow V$ is a diffeomorphism, then $(T_f\phi)(x)$ in $D'(V)$ is easy to define as
$|J_{f^{-1}}(x)|^{\frac{[\phi]}{d}}\phi(f^{-1}(x))$ in the sense of distributions as in~\cite{DuplantierRSV}.
Namely, for a test function $g\in D(V)$, define $(T_f\phi)(g)=\phi(T_{f^{-1}}g)$ where
$(T_{f^{-1}}g)(x)=|J_f(x)|^{1-\frac{[\phi]}{d}} g(f(x))$. However, this results in the unpleasant feature that domains $U$, $V$
will change with the map $f$. If one wants to implement a true action of the group $\cM(\RR^d)$
at the level of Schwartz distributions in $S'(\RR^d)$, one must stomach the point at infinity.
This was done in white noise theory (see~\cite{HidaKNY} and~\cite[\S 5.3]{Hida}) by changing the space of test functions
from $S(\RR^d)$ to the $D_{\chi}$ used in~\cite{GelfandGV} in order to revisit Bargmann's study~\cite{Bargmann}
of unitary representations of $\cM(\RR)$ and $\cM(\RR^2)$. A similar idea was used in a $p$-adic context in~\cite{LernerM,Lerner}.
However, most computations in physics are done in Fourier space and the ideal setting for Fourier analysis is provided by the
distribution space $S'(\RR^d)$. In order not to give up the beloved space $S'(\RR^d)$,
the following can be done when $0<[\phi]<\frac{d}{2}$.
If one does not second-guess Schwartz's intentions when choosing the letter ``S'' for his space, the explanation for the latter
is the word ``spherical''. Indeed, Schwartz proves (see~\cite{Schwartz1} or~\cite[\S VII.4]{Schwartz2}) that an element 
$\phi\in D'(\RR^d)$ belongs to $S'(\RR^d)$ if an only if it admits an extension to an element of $D'(\mathbb{S}^d)$ via
$\RR^d\subset\widehat{\RR^d}\simeq \mathbb{S}^d$
(see~\cite{AizenbudG} for a far-reaching generalization).
For a measure such as $\mu_{C_{-\infty}}$ or $\nu_{\phi}$, with $0<[\phi]<\frac{d}{2}$,
by the 1st-quantized KC theorem, the $\phi$'s would almost surely be in $\cC^{\alpha}$
for $\alpha<-[\phi]$. The key remark, is that by taking $\alpha>-d$, the extension is {\em unique}
(and it is also in the space $\mathcal{C}^{\alpha}$ for the sphere $\mathbb{S}^d$).
Namely, one is in a removable singularity situation similar to that of Riemann's Theorem (see~\cite[\S 8.1]{Dang}).
One can also realize the unique extension in a Borel measurable way by avoiding the Hahn-Banach Theorem
used in~\cite[\S VII.4]{Schwartz2} and relying instead on a more constructive
extension procedure (see~\cite[p. 88]{BogolubovLOT} or~\cite[Thm. 2.1]{Meyer}).
This gives a definition of the push-forward $(T_f)_{\ast}\nu$ of the original measure $\nu$
by an almost sure pathwise construction.
The additional restriction $[\phi]<\frac{d}{2}$, rather than just $[\phi]<d$, comes from the need for the DPC property in order
to reduce an equation $(T_f)_{\ast}\nu=\nu$ to a statement about pointwise correlations as in Conj. 4.
For the Gaussian measure $\mu_{C_{-\infty}}$, and contrary to what a hasty reading of~\cite{DuplantierRSV}
might suggest, CI (in the stronger
probabilistic sense just explained) holds for {\em all} $[\phi]\in (0,\frac{d}{2})$
because the DPC property
holds as well as the easy analogue of Conj. 4~\cite{Rajabpour}. The latter reduces to the fundamental
identity $|J(x)-J(y)|=|x-y|\times |x|^{-1}\times|y|^{-1}$ for the unit-sphere inversion.
It is somewhat surprising that CI needs the 1st-quantized KC Theorem and thus a condition such as $[\phi]<\frac{d}{2}$.
As for the ``second-quantized'' KC Theorem and the OPE, even for Gaussian fields, going beyond the
$[\phi]=\frac{d}{2}$ barrier is rather problematic.
This relates, e.g., to the construction of iterated integrals of fractional Brownian motion for low Hurst exponent.
Very little is known, apart from the findings from the difficult work by J. Unterberger, with help from
J. Magnen~\cite{Unterberger1,Unterberger2,MagnenU1,MagnenU2}.
Finally, before moving on to the $p$-adic fractional $\phi_{3}^{4}$ model, it is worth remarking that if CI and OS positivity hold, then so do more
exotic forms of the latter. One would have inversion positivity~\cite{FrankL,NeebO} with respect to say the unit sphere instead of a reflection
hyperplane.

\section{The $p$-adic hierarchical model}

It should be apparent from the previous discussion that Conj. 1-10 (except perhaps 5) are connected by a unifying theme: {\em robustness
with respect to space-dependent perturbations of the interaction}, e.g., replacing $g_r,\mu_r$ in (\ref{regularized}) by
functions $g_r(x),\mu_r(x)$.
The ultimate goal of this program is to make progress on these nine conjectures which is a battle that has to be fought {\em simultaneously}
on all fronts.
The chosen RG approach is very expensive: proofs such
as~\cite{BauerschmidtBS1,BauerschmidtBS2,BauerschmidtBS3,BauerschmidtBS4,BauerschmidtBS5,BrydgesS1,BrydgesS2,BrydgesS3,BrydgesS4,SladeT}
tend to take monumental proportions. Besides, the latter work addresses issues that are more tractable than the OPE,
anomalous dimensions and CI. The last thing one would want is to build a big RG machine and then, at the last
step towards proving say CI, to discover with horror that some choice of cut-off made at the beginning
will not work and one has to start over.
It is thus essential to find a toy model where {\em all} the problems in Conj. 1-10 are present
yet in {\em cleaner smaller-scale form}. This {\em is} what the $p$-adic hierarchical
fractional $\phi_{3}^{4}$ model of~\cite{AbdesselamCG1} provides.
When studying complex multiscale phenomena it is often important to decompose
functions into time-frequency atoms which live on a tree, e.g., when using a wavelet decomposition.
Unfortunately for most questions of interest the metric which governs how these atoms interact with each other is not
the natural (from the tree point of view) ultrametric distance, but the Euclidean metric of the {\em underlying continuum}.
Hierarchical models in physics amount to changing the model so it is the ultrametric distance which defines atomic correlations.
The same idea also appears in mathematics where such toy models are often called ``dyadic models''~\cite{Tao}.
Given a problem in Euclidean space, there are lots of ways of setting up a hierarchical model for it.
The $p$-adics, in some sense provide the most canonical, structured and principled way of doing this.
In order to make this article readable to probabilists who are not familiar with $p$-adic analysis, the $p$-adic fractional $\phi_{3}^{4}$ model will be introduced without appealing to $p$-adic analysis as covered, e.g.,
in~\cite{AlbeverioKS,Sally,VladimirovVZ} (see also~\cite{DragovichEtal} for a recent review) but instead by relying on the well developed probablistic intuition for trees (Galton-Watson
processes, branching Brownian motion, etc.). The few remarks needing $p$-adic analysis will be placed between the symbols $\spadesuit$
and $\heartsuit$.

Take an integer $p>1$ and let $\LL_k$ denote the set of boxes $\prod_{i=1}^{d}[p^k a_i, p^k(a_i+1))$, $(a_1,\ldots,a_d)\in\NN^d$, which form
a partition of the orthant $[0,\infty)^d$ by cubes of size $p^k$, $k\in\ZZ$.
The set $\cT=\cup_{k\in\ZZ}\LL_k$ naturally has the structure of a doubly-infinite tree. Now forget about how the tree
was defined using an artificial embedding in $[0,\infty)^d$ and just remember the {\em tree structure} and its organization
into {\em layers} $\LL_k$.
There is a clear notion of {\em rays} or {\em leaves at infinity} forming a set $\LL_{-\infty}=:\QQ_p^d$ (see~\cite{AbdesselamSlides}).
Pick a ray and call it the origin $0$. Give the name $\infty$ to the common root of all these rays.
For $x,y\in\QQ_p^d$ define their distance $|x-y|=p^k$
where $k$ labels the layer $\LL_k$ where the two rays split. This is an ultrametric, whose closed balls are in one-to-one
correspondence with vertices $\ux\in\cT$.
One can easily define a ``Lebesgue measure'' ${\rm d}^dx$ on $\QQ_p^d$ which assigns a volume $p^{kd}$ for a ball
corresponding to an $\ux$ in $\LL_k$. Thus, ${\rm Volume}=({\rm Linear\ size})^d$ as in $\RR^d$. The layer $\LL_0$
(the set of unit balls in $\QQ_p^d$) plays the role of the lattice $\ZZ^d$. Define the center $\Omega\in\cT$ as the site
where the bi-infinite path $\infty\rightarrow 0$ meets $\LL_0$.
Define a centered Gaussian random field $(\zeta_{\ux})_{\ux\in\cT}$
by letting $\EE\zeta_{\ux}\zeta_{\uy}$ if $\ux$, $\uy$ are in different layers or are in the same layer
$\LL_k$ but have different mothers. The remaining case is settled by letting $\EE\zeta_{\ux}\zeta_{\uy}$ be equal to $-p^{-d} p^{-k[\phi]}$
if $\ux\neq\uy$
and to $(1-p^{-d})p^{-k[\phi]}$ if $\ux=\uy$.
The Gaussian $\mu_{C_{-\infty}}$ is the law of the field $\phi(x)=\sum_{\ux} \zeta_{\ux}$ where the sum is over all sites along the ray
$\infty\rightarrow x$. It lives in $S'(\QQ_p^d)\simeq\RR^{\NN}$ obtained as the dual of the space
$S(\QQ_p^d)\simeq \mathfrak{s}_0\simeq \oplus_{\NN}\RR$ of {\em smooth} (i.e., {\em locally constant})
{\em compactly supported} functions $f:\QQ_p^d\rightarrow\RR$ with the locally convex topology defined by all seminorms.
Note that the usually trivial part of Bochner's Theorem is nontrivial but still true for $\mathfrak{s}_0$ since it is
sequential~\cite[pp. 795--796]{Nyikos}.
Fix $L=p^l$ for some integer $l>0$.
If one stops the summation over $\ux$, coming from $\infty$, at the layer $\LL_{lr}$, then this defines the analogue of the smooth field
$\phi\ast\rho_{{\rm UV},r}$ and its law $\mu_{C_r}$.
The volume cut-off $\rho_{{\rm IR},s}$ is the indicator function of the closed ball of radius $L^s$ around $0$.
Combine analogues of translations and $O(d)$ elements into the group of bijective isometries of $\QQ_p^d$ for the ultrametric.
Scaling correponds to vertical shifts between layers $\LL_k$.
Now the hierarchical analogues of Conj. 1-10 should be clear except for Conj. 3 and 4.
Define the cross-ratio {\em as before} for distinct points in $\widehat{\QQ_p^d}=\QQ_p^d\cup\{\infty\}$
{\em using the ultrametric} instead of the Euclidean metric. Define the group of global conformal transformations $\cM(\QQ_p^d)$
as the invariance group of the cross-ratio. Replace $|J_f(x)|$ by the Radon-Nikodym derivative relating the measure ${\rm d}^d x$
for $\QQ_p^d$ and its transform by $f$. Now Conj. 4 is clear. $\spadesuit$ Conj. 5 Part 1) degenerates since Fourier and lattice
cut-offs are the same but Part 2) remains. $\heartsuit$
Note that $\cT$ with the ``hyperbolic metric'' given by the {\em graph distance}
plays the role of $\mathbb{B}^{d+1}$ or $\mathbb{H}^{d+1}$.
The latter corresponds to the stratification of $\cT$ into layers $\LL_k$ whereas the conformal ball point of view involves stratification
into spheres (in graph distance) around the center $\Omega$.
Again (see~\cite[pp. 127--128]{LernerM}), {\em elements of $\cM(\QQ_p^d)$ correspond to hyperbolic isometries of $\cT$.}
The proof uses the beautiful Mumford-Manin-Drinfeld cross-ratio lemma (see~\cite[p. 246]{ManinD},~\cite[Lem. 5.6]{Manin1}
or~\cite[\S 3.2]{Manin2}):
${\rm CR}(x_1,x_2,x_3,x_4)=
p^{-\delta(x_{1}\rightarrow x_{2}; x_{3}\rightarrow x_{4})}$ where $\delta(x_{1}\rightarrow x_{2}; x_{3}\rightarrow x_{4})$ is the number
of common edges in the two paths, counted positively if orientations agree and negatively otherwise.
$\spadesuit$ The analogue of the inversion $J$ is $J(x)=|x|^{2} x$ with no typo! It preserves the cross-ratio
because of the fundamental identity $|J(x)-J(y)|=|x-y|\times |x|^{-1}\times|y|^{-1}$ which makes
the Euclidean norm special for $\RR^d$
and $|x|=\max_{1\le i\le d}|x_i|$ for $\QQ_p^d$.
$\heartsuit$
Finally, Conj. 3 just needs a definition of half-space and reflection. Take first $d=1$ and $p$ odd for simplicity
(one can also treat $p=2$). When drawing the tree $\cT$ and the path $\infty\rightarrow 0$ from bottom to top as in~\cite{AbdesselamSlides},
move half of the $p-1$ branches at each node to the left of the path and half to the right. This gives the analogue of
$\RR=(-\infty,0)\cup\{0\}\cup(0,\infty)$.
Then, take the product with $\QQ_p^{d-1}$
using the obvious interpretation of $\QQ_p^d$ as a cartesian product of $d$ copies of $\QQ_p^1$.
The resulting $p$-adic OS positivity is satisfied by $\mu_{C_{-\infty}}$ but not by $\mu_{C_r}$, just as in $\RR^d$.
Moreover, the range of fractional Laplacian exponents $\alpha$ for which OS positivity holds is the same
as that where $e^{-t(-\Delta)^\alpha}$
has a positive pointwise kernel~\cite{BlumenthalG}. One does {\em not}
have the limitation $\alpha\le 1$ on $\QQ_p^d$ which reflects
the simpler set of distributions supported at the origin and pole structure of the relevant Gamma functions.

\section{What has been done in the $p$-adic case}\label{oldwork}

For the $p$-adic model the following has been proved. The article~\cite{AbdesselamCG1} proves Conj. 1 except independence with respect to
$L$, $\rho_{\rm UV}$ and $\rho_{\rm IR}$ and full scale invariance. Only discrete scale invariance by powers of $L=p^{l}$
is proved.
Reference~\cite{AbdesselamCG1} also proves Conj. 6 with $[\phi^2]$ possibly depending on $L$ as well
as Conj. 7 except the DPC and CI
properties. Scale-invariance is also discrete.
Later in~\cite{AbdesselamCG2} (see~\cite[Ch. 4]{Chandra} for a preview), the possible dependence of $[\phi^2]$ on $L$ is {\em ruled out}
and discrete scale invariance
is {\em upgraded} to full one by powers of $p$ for Conj. 1 and Conj. 7.
In May 2015, the author proved the DPC property in Conj. 7.
Note that~\cite{AbdesselamCG1,AbdesselamCG2} not only prove conjectures but introduce {\em new methods}. The article~\cite{AbdesselamCG1}
develops the {\em space-dependent} rigorous RG, i.e., one able to handle couplings $g_r(x),\mu_r(x)$. The non-rigorous version
is called {\em local RG} in the physics literature (see~\cite[\S 12.4]{WilsonK} and~\cite{DrummondS,Osborn,JackO,Nakayama}).
The new method in~\cite{AbdesselamCG2}, based on a successful marriage of RG and correlation inequalities techniques,
shows that the measures $\nu_{\phi}$ obtained by the RG with two ratios $L$ and $p L$ are the same and thus inherit both discrete
scale invariances. This method should work in $\RR^d$ too since the ratio of $\log(L)$ and $\log(L+1)$ is irrational
(see, e.g.,~\cite{Waldschmidt}) and the DPC property would imply scale invariance by the full group $(0,\infty)$.

The first nontrivial example of successful use of the RG {\em method} is not due to Wilson but
rather Landen and Gauss in the late 18th century~\cite[\S 2.3]{McKeanM}. For $\vec{V}=(a,b)\in(0,\infty)^2$, let
$\cZ(\vec{V})=\int_{0}^{\pi/2} (a^2 \cos^2 \theta +b^2\sin^2 \theta)^{-\frac{1}{2}}\ {\rm d}\theta$ which is trivial to compute
on the ``line of fixed points'' $\{a=b\}$. The identity $\cZ(\vec{V})=\cZ(RG(\vec{V}))$ for the transformation
$RG:(a,b)\mapsto\left(\frac{a+b}{2},\sqrt{ab}\right)$ is the basis for Gauss's algorithm to compute $\cZ(\vec{V})$ in general.
The RG method in~\cite{AbdesselamCG1} (see~\cite{Abdesselam2013} for a quick review)
is based on the same idea.
Since the L\'evy Continuity Theorem holds for $\mathfrak{s}'_0\simeq\mathfrak{s}'_0\times\mathfrak{s}'_0$,
weak convergence follows from that of characteristic functions which is proved by showing uniform convergence
in a small {\em complex} neighborhood of the origin~\cite{Lukacs}.
The characteristic function is the ratio of $\cZ(\vec{V}[f,j])$ over  $\cZ(\vec{V}[0,0])=\cZ_{r,s}$ from
(\ref{regularized}) where $\vec{V}[0,0]$ is some data that encodes the integrand. Similarly, $\vec{V}[f,j]$
accounts for the addition of source terms $\phi(f)$ (or $\phi_{r,s}(f)$ in the notation of Conj. 6)
and $N_{r,s}[\phi^2](j)$ due to two test functions $f$ and $j$.
For suitable $\rho_{\rm UV}$, $\Gamma=C_0-C_1$ is positive semidefinite and one has a decomposition of Gaussian measures
$\mu_{C_0}=\mu_{\Gamma}\ast\mu_{C_1}$. The $RG$ transformation is the operation of performing the $\mu_{\Gamma}$ integration
followed by rescaling. It results in identities similar to that of Gauss. One then iterates, and carefully compares
the trajectories for the numerator and denominator using dynamical systems tools, some of which new like the
infinite-dimensional Poincar\'e-Koenigs Theorem of~\cite[\S9]{AbdesselamCG1}.
Although the author was unaware of it at the time when~\cite{AbdesselamCG1} was written, it turns out that~\cite[\S9]{AbdesselamCG1}
share a similar flavor with the treatment in~\cite[\S2.6 and \S2.7]{PalisDM} of the Stable Manifold Theorem and the $\lambda$-Lemma.
The main difference is that~\cite[\S2.6 and \S2.7]{PalisDM} is in the differentiable category whereas~\cite[\S9]{AbdesselamCG1}
is in the analytic category.
Other works in rigorous RG theory handle the $\phi(f)$ source by completing the square and translating the Gaussian measure.
This causes problems if one iterates $RG$ too many times. This is why~\cite{BrydgesDH,BauerschmidtBS1,Mitter} produce
scaling/continuum limits in finite volume. The new method of~\cite{AbdesselamCG1} allows one to take both $r\rightarrow -\infty$
{\em and} $s\rightarrow\infty$ limits and shows the
{\em commutation} of these limits. For the $N_{r,s}[\phi^2](j)$ source the previous translation
trick does not work and one has to introduce space-dependent RG techniques as in~\cite{AbdesselamCG1}. 
One could say that obtaining a CFT compactified on a torus~\cite{Mitter} is a good thing and
indeed it is~\cite[Ch. 10]{DiFrancescoMS},
provided one can extract the full space CFT from it by a short distance scaling limit~\cite[Eq. 3.2]{Dobrushin2}.
However, this exactly means doing $s\rightarrow\infty$ after $r\rightarrow -\infty$ which again seems to require techniques
as in~\cite{AbdesselamCG1}.

\section{What remains to be done}

In this last section we will try to isolate some specific tasks one needs to tackle in order to complete the program outlined in this article. Some must be done before others and it may thus be useful to indicate in which order one may proceed.
 
\noindent{\bf Phase 1:}
In this first phase the main task is to complete the proofs of the conjectures in the $p$-adic case.
Conj. 3 mostly needs the stable manifold theorem~\cite[\S 8.2]{AbdesselamCG1}.
Although technical, the OPE and the construction of pathwise squares form a natural continuation of~\cite{AbdesselamCG1}. One of the main issues to be addressed in this phase
concerns CI in Conj. 4 (CI in Conj. 7 follows from Conj. 4 and Conj. 8).
The claim in the truly remarkable article~\cite{Lerner} amounts to Conj. 4 but the article has a shortcoming: CI is shown for the ``wrong'' model. Namely, Lerner uses the {\em conformal ball model}
for his cut-off, i.e.,
spheres with center $\Omega$ instead of the layers $\LL_k$.
Our solution is simple: show that the ``right'' and ``wrong''
models are the same. We believe the main needed tools are in~\cite{AbdesselamCG1} whose space-dependent RG philosophy needs to be pushed one step
further by allowing {\em space-dependent cut-offs}!
Namely, one should reprove Conj. 1 with $r\rightarrow r(x)$ in (\ref{regularized})
and appropriate tweaks to the Lebesgue
measure involving $(\min\{1,|x|^{-1}\})^{2d}{\rm d}^d x$ which is the $p$-adic analogue of
$(1+|x|^2)^{-d}{\rm d}^d x$ defining the conformal Cauchy
distribution~\cite{Letac,DunauS}.

\noindent{\bf Phase 2:} 
One needs to tackle the {\em main problem}: developing a rigorous space-dependent RG for $\RR^d$ using work on the $p$-adic
side as a guide, following the philosophy in Weil's famous de profundis letter to his sister~\cite{Weil}.
The argument for CI in~\cite{PaulosRRZ} is in the sense of formal power series (FPS). As a warm-up, one should give a mathematical proof for this property (still in the FPS sense).
A thorough rigorous pertubative study of the fractional $\phi_{3}^{4}$ model on $\RR^3$ is available~\cite{BleherM1,BleherM2,Missarov1,Missarov2}
as well as an ingenious proof of CI on $\QQ_p$ in the FPS sense~\cite{LernerM}. One just needs to combine the two.
For the non-perturbative study, the first step is to find the right partition of unity for, e.g., $\int g(x) \phi(x)^4\ {\rm d}^d x$
in order to satisfy the {\em localization property} of~\cite{Abdesselam2013}. The sharp indicator functions used
in~\cite{BrydgesMS,AbdesselamCMP} do not look promising. One would need smooth decompositions, e.g., as in~\cite[Ch. 4]{AbdesselamX}.
The first thing to try are Daubechies wavelets which served Hairer well in ~\cite{Hairer}.
The wavelet transform will play a crucial role since it provides an extension to the upper half-space model.
A wavelet basis corresponds to sampling at the vertices of a $2$-adic $\cT$ embedded in the real $\mathbb{H}^{d+1}$.
One has to see how $\cM(\RR^d)$
acts on such a basis. There is a coupling/field duality in $\int g(x) \phi(x)^4\ {\rm d}^d x$ which may signal the need for different
dual treatments of $g$ and $\phi$, perhaps using~\cite{DaubechiesGM,FrazierJ}. 

\noindent{\bf Phase 3:} If enough progress is made on the main problem,
one may also explore possible applications to number theory.
By a remark of Burnol~\cite{Burnol1,Burnol2} and Bochner's Theorem for $D'((0,\infty))$~\cite[\S 6]{Fernique}, the Riemann Hypothesis (RH)
is {\em equivalent} to the existence of a
certain multiplicatively translation-invariant probability measure on $D'((0,\infty))$.
The covariance kernel behaves like $|x-y|^{-1}$ on the diagonal and the ideal KC Theorem would suggest the corresponding field has regularity
$\alpha<-\frac{1}{2}$ so one really needs random {\em distributions}.
While Burnol's explanation for Weil positivity is a quantum one (in the sense of quantum versus classical mechanics), that in 
the approaches by Connes~\cite{Connes} and Deninger~\cite{Deninger} are classical.
Let $G$ be a finite group with a left action on a finite set $X$ and denote by $g\mapsto T(g)$
the corresponding {\em unitary} representation in $l^2(X)$. Then, clearly, the matrix
$(M_{g_1,g_2})_{g_1,g_2\in G}$ given by $M_{g_1,g_2}={\rm tr}\left(T(g_1 g_2^{-1})\right)$ is positive
semidefinite. Connes' mechanism for positivity~\cite[p. 66]{Connes} is a sophisticated elaboration on the previous fact.
In Deninger's approach it is clear that RH is related to an action of the scaling group $(0,\infty)$, but action on what?
That is the question.
A possibility suggested in~\cite[p. 99]{Deninger} is that the relevant finite-dimensional dynamical system might emerge
as an attractor inside an infinite-dimensional one, similarly to renomalizable QFT's inside theory space where $(0,\infty)$ acts
as the RG. This tantalizing possibility was further investigated in~\cite{Leichtnam}.
Perhaps one should try to reconcile the quantum point of view of Burnol with the classical one of Connes and Deninger,
and holography~\cite{ManinM} may provide some clues for doing so. Although not yet very well understood, the holographic RG
is a way to geometrize the scale direction in Wilson's RG, roughly, using the $\{0\}^d\times (0,\infty)$ axis in $\mathbb{H}^{d+1}$
(see, e.g.,~\cite{deBoerVV,HeemskerkP,Nakayama}).
For instance, take a random field $\phi$ on $\RR^d$ and extend it to $\mathbb{H}^{d+1}$ using the wavelet transform.
This is the continuous analogue of the block-spinning procedure for the $p$-adic case in~\cite{LernerM,Lerner}
which goes back to~\cite{Zabrodin}.
Then restrict to the $(0,\infty)$ axis in $\mathbb{H}^{d+1}$. Multiplicative translation-invariance for the resulting
field, as in Burnol's remark, means that the original field is self-similar which most often than not~\cite{Nakayama}
translates into CI. However, this construction is too naive because the extended field is smooth which is incompatible
with the expected $-\frac{1}{2}$ H\"{o}lder regularity. In any case, a good example to look at would be the simplest nontrivial
case of elliptic curves over finite fields where RH is not a conjecture but a theorem from 1936~\cite{Hasse}.
Quite ironically, the profound unity in the book series by Gel'fand and co-authors on generalized functions
was {\em literally} lost in translation.
Volumes 5~\cite{GelfandGV} and 6~\cite{Gelfand6}
focused on the unitary representations of groups such as $\cM(\RR)$ or $\cM(\RR^2)$, over local fields such as $\RR$ or $\QQ_p$.
The corresponding applications to number theory form a flourishing and vibrant area in mathematics (see~\cite{Frenkel} for
review emphasizing connections to CFT).
It would be desirable to develop number theoretical applications, preferably easier ones
than the RH bogeyman, for Volume 4~\cite{Gelfand4} too.

\noindent{\bf Acknowledgements:}
The author would like to thank, first and foremost, his collaborators A. Chandra and G. Guadagni without whom~\cite{AbdesselamCG1,AbdesselamCG2}
would not have taken off the ground.
For
useful face-to-face discussions, email correspondence
or posts on MathOverflow, the author thanks: D. Buchholtz,
R. Budney, P. Cartier, N. V. Dang, M. Emerton, M. Furlan, M. Hairer, S. Krushkal,
P. Michor, J.-C. Mourrat, R. Rhodes, S. Rychkov, V. Vargas and O. Zeitouni.
The work~\cite{AbdesselamCG1} reported on in this article was supported in part by the National Science Foundation under grant DMS \# 0907198.
Last but not least, the author warmly thanks the organizers of the 6th International Conference on $p$-adic Mathematical Physics and its Applications, and in particular W. Zuniga-Galindo, for a very stimulating workshop.


\end{document}